\documentclass[oneside,english,reqno]{amsart}
\usepackage[T1]{fontenc}
\usepackage[latin9]{inputenc}
\usepackage{mathrsfs}
\usepackage{amsbsy}
\usepackage{amstext}
\usepackage{amsthm}
\usepackage{amssymb}
\usepackage{stackrel}
\usepackage{esint}

\makeatletter
\numberwithin{equation}{section}
\numberwithin{figure}{section}
\theoremstyle{remark}
\newtheorem*{rem*}{\protect\remarkname}
\theoremstyle{plain}
\newtheorem*{thm*}{\protect\theoremname}
\theoremstyle{plain}
\newtheorem{thm}{\protect\theoremname}
\theoremstyle{plain}
\newtheorem{lem}[thm]{\protect\lemmaname}

\usepackage{etoolbox}

\patchcmd{\@maketitle}
  {\ifx\@empty\@dedicatory}
  {\smallskip
    \begin{center}
    \footnotesize
    \begin{tabular}{c}
      16110180025@fudan.edu.cn \\
      Department of Mathematics \\
      Fudan University \\
      Shanghai, China \\
	  June, 2018
    \end{tabular}
  \end{center}
  \ifx\@empty\@dedicatory}
  {}{}

\makeatother

\usepackage{babel}
\providecommand{\lemmaname}{Lemma}
\providecommand{\remarkname}{Remark}
\providecommand{\theoremname}{Theorem}

\begin{document}
\title[Limit of 1D Mixed-Mechanism IPS Model]{Limit of 1-Dimensional Mixed-Mechanism Interacting Particle System
Model }
\author{Tong Zhao}
\maketitle

\part*{Introduction}

Elaborating on the model from voter process with mixed-mechanism under
suitable scaling, I have two new mechanisms which are random switch
and unbiased local Homogenization and subtly biased advantage but
with state dependent coefficient involved. The most crucial one, the
existence of high-frequency duplication generating the diffusion term
and noise term in each case identifies the limit equation as SPDE
driven by space time white noise. 

At the beginning, presenting the notations (parameters in the model)
used in this paper is necessary for convenience. $\mathcal{N}_{n}$
is the number of neighbors of a voter, $\mathcal{D}_{n}$ is the distance
between two voters who neighbor each other at most, $\mathcal{H}_{n}$
is the rate of high-frequency mechanism, $\mathcal{L}_{n}$ is the
rate of low-frequency mechanism, $\mathcal{S}_{n}$ is the scale of
density and $\mathcal{\rho}_{n}:=\mathcal{N}_{n}/(2\mathcal{D}_{n})$
is the density of voters in \textit{the n}th mode. 

For any fixed $n\in\mathbb{Z}$, a classic 1-dimensional model is
the lattice scale is determined by $\mathbb{Z}/\mathcal{\rho}_{n}$
in its \textit{n}th model. There is a voter on each lattice who is
an advocate of either $\boldsymbol{\mathbf{A}}\mathit{}$ or $\boldsymbol{a}$
. If $x,y\in\mathbb{Z}/\mathcal{\rho}_{n}$ and $|x-y|\leqslant\mathcal{D}_{n}$,
we regard $x,y$ as neighbors denoted by $x\sim y$. 

We will use $d_{n}(t,x)$, $u_{n}(t,x)$ and $u(t,x)$ to denote the
density, the $\mathcal{S}_{n}$-scaled density in the \textit{n}th
model and its limit of advocates of $\boldsymbol{\mathbf{A}}$ at
$x$ point and time $t$respectively. We can also consider this model
on a ring $(\mathbb{Z}/n)/\mathbb{Z}$ as its \textit{n}th model.
Under this circumstance, we need $u$ and $\dot{W}$ of period $1$
in final equation.

In a particular model: $\mathcal{N}_{n}=2n^{1/2},\mathcal{D}_{n}=n^{-1/2},\mathcal{H}_{n}=2n,\mathcal{L}_{n}=2\theta,\mathcal{S}_{n}=1,\mathcal{\rho}_{n}=n,d_{n}=\mathcal{S}_{n}u_{n}$.

\part*{Mechanisms}

The initiative, 1-dimensional voter process can converge to a SPDE
driven by space time white noise with various drift terms---especially
the bistable drift term--- in our model in vague sense, was motivated
by Allee effect. In time, we find it can be generalized to voter process
and then go a step further to a more general form. 

\section{High-Frequency Unbiased Oscillation Mechanism}

\subsection{Symmetric Duplication}

The reason we call it high-frequency unbiased oscillation mechanism
is its rate must be high enough to generate the Laplacian term and
the white noise term and is not related to the number of neighbors
of each individual on $\mathbb{Z}/n$. For anyone of two individuals
$x,y$ neighboring one another adopting the view of the other in a
period of time (i.e. $[0,t]$) indefinitely times subject to Poisson
process at rate $\mathcal{H}_{n}/\mathcal{N}_{n}$ in \textit{the
n}th model, these Poisson processes labeled by $P_{t}(x,y)$ are independent
of each other with ordered pairs $(x,y)$. 

Let consider how this mechanism gives birth to laplacian term and
white noise term. We define:
\[
\Delta_{n}(f)(x):=\mathcal{H}_{n}\sum_{y\sim x}\frac{1}{\mathcal{N}_{n}}(f(y)-f(x)).
\]

Whence $\text{\ensuremath{\Delta}}_{n}$ is a generator of a random
walk at rate $\mathcal{H}_{n}$ , which is the totally rate of diffusion
of density if substituting approximate local density $u_{n}$ for
$f$. The transition probability of this random walk is uniformly
distributed on the lattices in the neighborhood of a specific point.
So the variance of the distribution of the transition probability
is $\mathcal{D}_{n}^{2}/3$. Since each jump of the random walk is
independent of one another and the expectation of the number of jumps
in unit time is $\mathcal{H}_{n}$, the variance of the random walk
in unit time is $\mathcal{D}_{n}^{2}\cdot\mathcal{H}_{n}/3$. Therefore,
we require the convergence of $\mathcal{D}_{n}^{2}\cdot\mathcal{H}_{n}/3$
upon $n$ tending to infinity, assuming this limit is $\sigma^{2}$.
According to local central limit theorem, we have the laplacian term
$\frac{\sigma^{2}}{2}\text{\ensuremath{\Delta}}u$ . 

As to the white noise term, so $d_{n}\rho_{n}$ is the density of
$\boldsymbol{\mathbf{A}}$ in the neighborhood of a certain point
x (i.e. the number of $\boldsymbol{\mathbf{A}}$ on a unit interval
of space). In a unit interval of time, the expectation of the number
of times of occurrence of adopting view of others on a unit interval
whose density is uniform and identified with $\mathcal{\rho}_{n}$
is $\mathcal{H}_{n}\mathcal{\rho}_{n}$. The distribution of increment
of the density upon the next adopting occurring is: 

\[
\begin{array}{ccc}
-1 & 0 & 1\\
(1-d_{n})d_{n} & 1-2d_{n}(1-d_{n}) & d_{n}(1-d_{n})
\end{array}
\]

Hence we have:
\[
\partial_{t}(d_{n}\mathcal{\mathcal{\rho}}_{n})=\sqrt{2d_{n}(1-d_{n})(\mathcal{H}_{n}\mathcal{\rho}_{n})}\dot{W},
\]

which leads to 
\[
\partial_{t}u_{n}=\sqrt{2u_{n}(1-u_{n}\mathcal{S}_{n})\mathcal{H}_{n}/(\mathcal{S}_{n}\mathcal{\rho}_{n})}\dot{W}.
\]

There we require some convergence property to get a non-trivial term. 

\section{Low-Frequency Drift Mechanism}

Now, we will pay attention to the significant part, low-frequency
drift mechanism. The following are various types of them, whose occurrence
related to a lattice point $x$ is of the same order as $\mathcal{L}_{n}$.

\subsection{Mutation (random switching)}

This mechanism in voter process model means every voter switches their
views or changes their minds somehow, maybe on a whim or something
irrational. Assuming the rate of the mutation of a kind of voter,
$\mathcal{L}_{n}$, we are handy to conclude that if $\boldsymbol{\mathbf{A}}$
mutates into $\boldsymbol{a}$at rate $\mathcal{L}_{n}$, we obtain
the term $-u_{n}\mathcal{L}_{n}$, and in turn if $\boldsymbol{a}$mutates
into $\boldsymbol{\mathbf{A}}$ at rate $\mathcal{L}_{n}$, we obtain
the term $(\mathcal{S}_{n}^{-1}-u_{n})\mathcal{L}_{n}$ by observing
the quantity of voter $\boldsymbol{\mathbf{A}}$ in density sense.
Specially, if the mutation is unbiased, we have a term $(\mathcal{S}_{n}^{-1}-2u_{n})\mathcal{L}_{n}$
.

\subsection{Asymmetric Duplication}

We always come into a situation where advocates of the party, $\boldsymbol{\mathbf{A}}$,
are more active and aggressive than the other $\boldsymbol{a}$. That
means $\boldsymbol{\mathbf{A}}$ has an extra frequency to change
$\boldsymbol{a}$'s mind or $\boldsymbol{\mathbf{A}}$'s view is more
easily adopted by $\boldsymbol{a}$. Similar to high-frequency mechanism,
we use $\widetilde{P}_{t}(x,y)$ to represent the process that account
the number of times of $x$trying adopting $y$'s opinion in time
interval $[0,t]$, but only if y is an advocate of $\boldsymbol{\mathbf{A}}$
when $\widetilde{P}_{t}(x,y)$ jumps, x succeeds in adopting y's opinion
$\boldsymbol{\mathbf{A}}$. $\widetilde{P}_{t}(x,y)$ is subject to
Poisson process with expectation $\mathcal{L}_{n}/\mathcal{N}_{n}$.
Same as the above considering a neighborhood of a point x, we have:

\[
\partial_{t}(d_{n}\mathcal{\rho}_{n})=(d_{n}\mathcal{\rho}_{n})\mathcal{L}_{n}/\mathcal{N}_{n}\cdot\mathcal{N}_{n}(1-d_{n}).
\]

Hence,

\[
\partial_{t}u_{n}=\mathcal{L}_{n}u_{n}(1-u_{n}\mathcal{S}_{n}).
\]

\begin{rem*}
This mechanism could not be symmetric, otherwise
\end{rem*}
\[
\partial_{t}u_{n}=\mathcal{L}_{n}u_{n}(1-u_{n}\mathcal{S}_{n})-\mathcal{L}_{n}(1-u_{n}\mathcal{S}_{n})u_{n}=0.
\]

This implies compensation, hence no effect.

\subsection{Local Homogenization ( multi-consulting )}

Some voters may be very stubborn and discreet or even scrupulous,
their points of views are more stable. For example, every stubborn
individual $x$ inquiring two neighbors $y,z$ several times subject
to Poisson process with a fixed rate $\mathcal{L}_{n}/\mathcal{N}_{n}^{2}=\widetilde{P}_{t}(x,y,z)$
about their views in the time interval $\left[0,t\right]$ will alter
his or hers at a jumping moment of $\widetilde{P}_{t}(x,y,z)$ only
if theirs all differ from his or hers at that time. If voters for
$\boldsymbol{\mathbf{A}}$ are stubborn, we considering the neighborhood
of a fixed point on lattice obtain:

\[
\partial_{t}(d_{n}\mathcal{\rho}_{n})=(-d_{n}\mathcal{\rho}_{n})\mathcal{L}_{n}/\mathcal{N}_{n}^{2}\cdot\mathcal{N}_{n}(1-d_{n})\cdot\mathcal{N}_{n}(1-d_{n}).
\]

Thereby:

\[
\partial_{t}u_{n}=-\mathcal{L}_{n}u_{n}(1-u_{n}\mathcal{S}_{n})^{2}.
\]

However, if voters for $\boldsymbol{a}$ are also stubborn, we obtain 

\[
\partial_{t}u_{n}=\mathcal{L}_{n}\left[((1-u_{n})u_{n}^{2}-u_{n}(1-u_{n})^{2}\right]=2\mathcal{L}_{n}u_{n}(u_{n}-\frac{1}{2})(1-u_{n})
\]

similarly without loss of the highest order term, which leads to a
symmetric bistable structure if $\mathcal{L}_{n}\rightarrow constant,\mathcal{S}_{n}\rightarrow1$.

\subsection{Polynomial }

Following the above suit, we consider the case every stubborn individual
$x$ inquiring $m$ neighbors $y_{1},y_{2},\cdots,y_{m}$ several
times subject to Poisson process with a fixed rate $\mathcal{L}_{n}/\mathcal{N}_{n}^{m}=\widetilde{P}_{t}(x,y_{1},y_{2},\cdots,y_{m})$
about their views in the time interval $\left[0,t\right]$ will alter
his or hers at a jumping moment of $\widetilde{P}_{t}(x,y_{1},y_{2},\cdots,y_{m})$
only if theirs all differ from his or hers at that time. Similarly,
we obtain an m-order polynomial in drift term.

If $\mathcal{L}_{n}\rightarrow constant,\mathcal{S}_{n}\rightarrow1$:
\begin{enumerate}
\item As in the above case, we know $\frac{1}{2}$ is the only zero point
between 0 and 1 whatever $m$ is, which also leads to a symmetric
bistable structure.
\item and if we change m and modify the relative intensity to make $\boldsymbol{\mathbf{A}}$
and $\boldsymbol{a}$ is not symmetric, we can have a general formula:
\[
\partial_{t}u_{n}=P(u_{n}),
\]
where $P(\cdot)$ is a polynomial. 
\end{enumerate}

\subsection{State dependence }

In some cases, the rate $\mathcal{L}_{n}$ of a point $x$ may depend
on the state of its neighborhood. Considering $\mathcal{L}_{n}(u_{n})$
as a function of $u_{n}$, we apply this assumption to 2.1, 2.2 (others
can be prove in the same way)and capture the equations below: 

\[
\partial_{t}u_{n}\mathcal{S}_{n}=-u_{n}\mathcal{S}_{n}\mathcal{L}_{n}(u_{n})
\]

\[
\partial_{t}(d_{n}\mathcal{\rho}_{n})=(d_{n}\mathcal{\rho}_{n})\mathcal{L}_{n}(u_{n})/\mathcal{N}_{n}\cdot\mathcal{N}_{n}(1-d_{n}),
\]
i.e.

\[
\partial_{t}u_{n}=-u_{n}\mathcal{L}_{n}(u_{n})
\]

\[
\partial_{t}u_{n}=\mathcal{L}_{n}(u_{n})u_{n}(1-u_{n}\mathcal{S}_{n}).
\]
 If $\mathcal{S}_{n}\rightarrow constant,\mathcal{L}_{n}(u_{n})=P(u_{n})$
:

\[
\partial_{t}u=-uP(u)
\]

\[
\partial_{t}u=P(u)u(1-u).
\]

\part*{Preliminaries and Description of the Theorem}

Choosing different mechanism and modifying parameters at the beginning
appropriately to warrant the convergence, we would get a non-trivial
SPDE.

Set :
\begin{enumerate}
\item $D(f,\delta)(x):=sup\left\{ \left|f\left(y\right)-f\left(x\right)\right|:\left|y-x\right|\leqslant\delta,y\in\mathbb{Z}/\mathbb{\rho}_{n}\right\} $,
for $x\in\mathbb{Z}/\mathbb{\rho}_{n},\delta>0$.
\item $\xi_{t}^{n}(x)\in\left\{ 0,1\right\} $ is an identifier of state
of voter at $x\in\mathbb{Z}/\mathbb{\rho}_{n}$ and time $t$, without
lose of generality, $\boldsymbol{\mathbf{A}}$, $\boldsymbol{a}$
corresponding to 1 and 0 respectively. Then the dynamics of $\xi_{t}^{n}(x)$
are noted according to various mechanism.
\item $A(f)(x):=\left(\mathcal{N}_{n}\mathcal{S}_{n}\right)^{-1}\underset{y\sim x}{\sum}f(y)$
for $x\in\mathbb{Z}/\mathbb{\rho}_{n}$ , then linearly interpolated.
\item $u_{n}(t,x):=A(\xi_{t}^{n})(x)$ .
\item $\hat{u}_{n}(t,x):=u_{t}(x)-(v_{0},\psi_{t}^{x}).$
\item $\tilde{u}_{n}(t,x):=\hat{u}_{n}(t,x)$ on the lattice $z\in\mathbb{Z}/\mathcal{\rho}_{n},\;t\in\mathbb{N}/(n\mathcal{\rho}_{n})$,
then linearly interpolate first in $x$ and then in $t$ to obtain
a continuous $\mathcal{\mathscr{C}}$ valued process.
\item $v_{t}^{n}(x):=$$\left(\mathcal{\rho}_{n}\mathcal{S}_{n}\right)^{-1}\underset{x}{\sum}\delta_{x}I\text{\ensuremath{\left(\xi_{t}^{n}(x)=1\right)}}$
the measure valued process.
\item $(f,g):=\mathcal{\rho}_{n}^{-1}\underset{x}{\sum}f(x)g(x)$ for $f,g:\mathbb{Z}/\mathcal{\rho}_{n}\longrightarrow\mathfrak{\mathbb{R}}$
and $(v,f)=\int fdv$ for $f:\mathbb{Z}/\mathcal{\rho}_{n}\longrightarrow\mathfrak{\mathbb{R}},v\in\mathcal{M}(\mathbb{Z}/\mathcal{\rho}_{n})$.
\item $e_{\lambda}(x):=e^{\lambda\left|x\right|}$ for all $x\in\mathcal{\mathbb{R}}$.
\item $\mathscr{C:=}\left\{ f:\mathbb{R\rightarrow}[0,\infty)\;continuous\;with\;\left|f(x)e_{\lambda}(x)\right|\rightarrow0\;as\;\left|x\right|\rightarrow\infty\;for\;all\;\lambda<0\right\} $
.
\item $\left\Vert f\right\Vert _{\lambda}:=\underset{x}{\sup}\left|f(x)e_{\lambda}(x)\right|$
\item $\underset{p(x,y,z)}{\sum}f(x,y,z):=f(x,y,z)+f(z,x,y)+f(y,z,x)$
\item $p(\sigma^{2},\cdot)$ is density function of centered normal distribution
with variance $\sigma^{2}$.
\item $C$ , a constant having nothing to do with our interest, is different
from line to line.
\item $\left|P(\cdot)\right|:=\stackrel[i=0]{n}{\sum}\left|p_{i}\right|,\quad if\;P(x)=\stackrel[i=0]{n}{\sum}p_{i}x^{i}.$
\end{enumerate}
In the following argument, we will omit superscript $n$ without ambiguity.
\begin{thm*}
Upon $n$ tend to infinity, and $u_{n}(0,x)$ converge to $f_{0}$
in $\mathscr{C}$ sense. Then $u_{n}(t,x)$ converge in distribution
sense to a continuous $\mathscr{C}$ valued process $u(t,x)$ which
solves the following SPDEs under respective conditions.
\end{thm*}
\begin{enumerate}
\item If you choose symmetric duplication as high-frequency and multiple
consulting and asymmetric state dependent mutation as low-frequency
mechanism respectively, then tune parameters to $\mathcal{N}_{n}=2n^{1/2},\mathcal{D}_{n}=n^{-1/2},\mathcal{H}_{n}=2n,\mathcal{L}_{n}^{c}=2,\mathcal{L}_{n}^{m}=P(u_{n}),\mathcal{S}_{n}=1,\mathcal{\rho}_{n}=n,d_{n}=\mathcal{S}_{n}u_{n},$
$u(t,x)$ will fit 
\[
\partial_{t}u=\frac{1}{3}\text{\ensuremath{\Delta}}u+4u(u-\frac{1}{2})(1-u)-uP(u)+2\sqrt{u(1-u)}\dot{W},\;u(0,x)=f_{0}.
\]
\item If you choose symmetric duplication and asymmetric state dependent
mutation as high-frequency and low-frequency mechanism respectively,
then tune parameters to $\mathcal{N}_{n}=2n^{1/2},\mathcal{D}_{n}=n^{-1/2},\mathcal{H}_{n}=2n,\mathcal{L}_{n}(u_{n})=P(u_{n}),\mathcal{S}_{n}=1,\mathcal{\rho}_{n}=n,d_{n}=\mathcal{S}_{n}u_{n}$,
$u(t,x)$will fit
\[
\partial_{t}u=\frac{1}{3}\text{\ensuremath{\Delta}}u-uP(u)+2\sqrt{u(1-u)}\dot{W},\qquad u(0,x)=f_{0}.
\]
\item If you choose asymmetric duplication and mutation as high-frequency
mechanism and no low-frequency mechanism , then tune parameters to
$\mathcal{N}_{n}=2n^{3/2},\mathcal{D}_{n}=n^{-1/2},\mathcal{H}_{n}=n,\mathcal{L}_{n}=0,\mathcal{S}_{n}=n^{-1},\mathcal{\rho}_{n}=n^{2},d_{n}=\mathcal{S}_{n}u_{n}$,
$u(t,x)$ will fit
\[
\partial_{t}u=\frac{1}{6}\text{\ensuremath{\Delta}}u+\sqrt{2u}\dot{W},\qquad u(0,x)=f_{0}.
\]
\item If you only choose symmetric mutation as high-frequency and no low-frequency
mechanism, then tune parameters to $\mathcal{N}_{n}=2n^{1/2},\mathcal{D}_{n}=n^{-1/2},\mathcal{H}_{n}=n,\mathcal{L}_{n}=0,\mathcal{S}_{n}=1,\mathcal{\rho}_{n}=n,d_{n}=\mathcal{S}_{n}u_{n}$,
$u(t,x)$ will fit
\[
\partial_{t}u=\dot{W,}\qquad u(0,x)=f_{0}=\frac{1}{2}.
\]
 The conditions and parameters are representative, you can follow
my suit to derive similar SPDE with various combinations of mechanisms
from following process. Therefore we will only give the proof of case
1.
\end{enumerate}
All mechanisms of one kind are subject to i.i.d Poisson processes,
and process between different kinds are independent mutually:

voter $x$ takes value of $y$

\[
(P_{s}(x,y):x,y\in\mathbb{Z}/\mathcal{\rho}_{n},x\sim y)\quad with\;rate\;\mathcal{H}_{n}/\mathcal{N}_{n}=n^{1/2},
\]

state dependent mutation from $\boldsymbol{\mathbf{A}}\mathit{}$
to $\boldsymbol{a}$ 

\[
(\hat{P}_{s}(x):x\in\mathbb{Z}/\mathcal{\rho}_{n})\quad with\;rate\;P(u_{n}),
\]

voter x consults $y$and $z$

\[
(\tilde{P}_{s}(x,y,z):x,y,z\in\mathbb{Z}/\mathcal{\rho}_{n},y\sim x,z\sim x)\quad with\;rate\;\mathcal{L}_{n}^{c}/\mathcal{N}_{n}^{2}=2\theta n^{-1}.
\]

The dynamics of the process in case 1 is described below:

\begin{align*}
\xi_{t}^{n}(x) & =\xi_{0}^{n}(x)+\sum_{y\sim x}\intop_{0}^{t}(\xi_{s-}^{n}(y)-\xi_{s-}^{n}(x))dP_{s}(x,y)-\intop_{0}^{t}\xi_{s-}^{n}(x)d\hat{P}_{s}(x)\\
 & +\sum_{y\sim x}\sum_{z\sim x}\intop_{0}^{t}(1-\xi_{s-}^{n}(x))\xi_{s-}^{n}(y)\xi_{s-}^{n}(z)d\tilde{P}_{s}(x,y,z)\\
 & -\sum_{y\sim x}\sum_{z\sim x}\intop_{0}^{t}\xi_{s-}^{n}(x)(1-\xi_{s-}^{n}(y))(1-\xi_{s-}^{n}(z))d\tilde{P}_{s}(x,y,z)\\
 & =\xi_{0}^{n}(x)+\sum_{y\sim x}\intop_{0}^{t}(\xi_{s-}^{n}(y)-\xi_{s-}^{n}(x))dP_{s}(x,y)-\intop_{0}^{t}\xi_{s-}^{n}(x)d\hat{P}_{s}(x)\\
 & -\sum_{y\sim x}\sum_{z\sim x}\intop_{0}^{t}\xi_{s-}^{n}(x)d\tilde{P}_{s}(x,y,z)\\
 & +\sum_{y\sim x}\sum_{z\sim x}\intop_{0}^{t}(\sum_{p(x,y,z)}\xi_{s-}^{n}(x)\xi_{s-}^{n}(y))d\tilde{P}_{s}(x,y,z)\\
 & -2\sum_{y\sim x}\sum_{z\sim x}\intop_{0}^{t}\xi_{s-}^{n}(x)\xi_{s-}^{n}(y)\xi_{s-}^{n}(z)d\tilde{P}_{s}(x,y,z).
\end{align*}

Then take a test function $\phi:[0,\infty)\times\mathbb{Z}/\mathcal{\rho}_{n}\rightarrow\mathbb{R}$
with $t\rightarrow\phi_{t}(x)$ continuously differentiable and satisfying
\[
\intop_{0}^{T}\left(\left|\phi_{s}\right|+\phi_{s}^{2}+\left|\partial_{s}\phi_{s}\right|,1\right)ds<\infty.
\]

Implementing integration by parts, for $t\leqslant T$, we have
\begin{align*}
(v_{t},\phi_{t}) & =(v_{0},\phi_{0})+\intop_{0}^{t}\left(v_{s},\partial_{s}\phi_{s}\right)ds\\
 & +\left(\mathcal{\rho}_{n}\mathcal{S}_{n}\right)^{-1}\sum_{x}\sum_{y\sim x}\intop_{0}^{t}\xi_{s-}^{n}(y)(\phi_{s}(x)-\phi_{s}(y))dP_{s}(x,y)\\
 & -\left(\mathcal{\rho}_{n}\mathcal{S}_{n}\right)^{-1}\sum_{x}\intop_{0}^{t}\xi_{s-}^{n}(x)\phi_{s}(x)d\hat{P}_{s}(x)\\
 & -\left(\mathcal{\rho}_{n}\mathcal{S}_{n}\right)^{-1}\sum_{x}\sum_{y\sim x}\sum_{z\sim x}\intop_{0}^{t}\xi_{s-}^{n}(x)\phi_{s}(x)d\tilde{P}_{s}(x,y,z)\\
 & +\left(\mathcal{\rho}_{n}\mathcal{S}_{n}\right)^{-1}\sum_{x}\sum_{y\sim x}\sum_{z\sim x}\intop_{0}^{t}(\sum_{p(x,y,z)}\xi_{s-}^{n}(x)\xi_{s-}^{n}(y))\phi_{s}(x)d\tilde{P}_{s}(x,y,z)\\
 & -2\left(\mathcal{\rho}_{n}\mathcal{S}_{n}\right)^{-1}\sum_{x}\sum_{y\sim x}\sum_{z\sim x}\intop_{0}^{t}\xi_{s-}^{n}(x)\xi_{s-}^{n}(y)\xi_{s-}^{n}(z)\phi_{s}(x)d\tilde{P}_{s}(x,y,z)\\
 & +\left(\mathcal{\rho}_{n}\mathcal{S}_{n}\right)^{-1}\sum_{x}\sum_{y\sim x}\intop_{0}^{t}\xi_{s-}^{n}(x)\phi_{s}(x)(dP_{s}(y,x)-dP_{s}(x,y)).
\end{align*}

\section{Semi-Martingale Decomposition}

\subsection{Laplacian Term}

We break the term into two parts, a fluctuation term and an average
term

\begin{align*}
 & \left(\mathcal{\rho}_{n}\mathcal{S}_{n}\right)^{-1}\sum_{x}\sum_{y\sim x}\intop_{0}^{t}\xi_{s-}^{n}(y)(\phi_{s}(x)-\phi_{s}(y))(dP_{s}(x,y)-(\mathcal{H}_{n}/\mathcal{N}_{n})ds)\\
 & +\left(\mathcal{\rho}_{n}\mathcal{S}_{n}\right)^{-1}\sum_{x}\sum_{y\sim x}\intop_{0}^{t}\xi_{s-}^{n}(y)(\phi_{s}(x)-\phi_{s}(y))(\mathcal{H}_{n}/\mathcal{N}_{n})ds\\
 & =E_{t}^{(1)}(\phi)+\intop_{0}^{t}(v_{s},\Delta_{n}(\phi_{s}))ds,
\end{align*}

where
\begin{align*}
E_{t}^{(1)}(\phi) & :=\left(\mathcal{\rho}_{n}\mathcal{S}_{n}\right)^{-1}\sum_{x}\sum_{y\sim x}\intop_{0}^{t}\xi_{s-}^{n}(y)(\phi_{s}(x)-\phi_{s}(y))(dP_{s}(x,y)-(\mathcal{H}_{n}/\mathcal{N}_{n})ds)
\end{align*}

is a martingale with brackets process given by

\begin{align*}
d\left\langle E_{n}^{(1)}(\phi)\right\rangle _{t} & =\left(\mathcal{\rho}_{n}\mathcal{S}_{n}\right)^{-2}\sum_{x}\sum_{y\sim x}\xi_{t}^{n}(y)(\phi_{t}(x)-\phi_{t}(y))^{2}(\mathcal{H}_{n}/\mathcal{N}_{n})dt\\
 & \leqslant\left(\mathcal{\rho}_{n}\mathcal{S}_{n}\right)^{-2}\sum_{y}\xi_{t}^{n}(y)D^{2}(\phi_{t},\mathcal{D}_{n})(y)\mathcal{H}_{n}dt\\
 & =\left(\mathcal{\rho}_{n}\mathcal{S}_{n}\right)^{-1}\mathcal{H}_{n}(v_{t},D^{2}(\phi_{t},\mathcal{D}_{n}))dt\\
 & \leqslant\left(\mathcal{\rho}_{n}\mathcal{S}_{n}\right)^{-1}\mathcal{H}_{n}\left\Vert D(\phi_{t},\mathcal{D}_{n})\right\Vert _{\lambda}^{2}(e_{-2\lambda},v_{t})dt\\
 & \leqslant\left(\mathcal{\rho}_{n}\mathcal{S}_{n}\right)^{-1}\mathcal{H}_{n}\left\Vert D(\phi_{t},\mathcal{D}_{n})\right\Vert _{\lambda}^{2}(e_{-2\lambda},1)dt.
\end{align*}

Alternatively, we bound it by

\begin{align*}
d\left\langle E_{n}^{(1)}(\phi)\right\rangle _{t} & \leqslant\left(\mathcal{\rho}_{n}\mathcal{S}_{n}\right)^{-2}\sum_{x}\sum_{y\sim x}2\xi_{t}^{n}(y)\left\Vert \phi_{t}\right\Vert _{0}(\left|\phi_{t}(x)\right|+\left|\phi_{t}(y)\right|)(\mathcal{H}_{n}/\mathcal{N}_{n})dt\\
 & =2\left\Vert \phi_{t}\right\Vert _{0}\left(\mathcal{\rho}_{n}\mathcal{S}_{n}\right)^{-1}\mathcal{H}_{n}\left[(\left|\phi_{t}\right|,u_{t}\mathcal{S}_{n})+(v_{t},\left|\phi_{t}\right|)\right]dt\\
 & \leqslant4\left\Vert \phi_{t}\right\Vert _{0}\left(\mathcal{\rho}_{n}\mathcal{S}_{n}\right)^{-1}\mathcal{H}_{n}(\left|\phi_{t}(x)\right|,1)dt.
\end{align*}

\subsection{Drift Term}

We break the term into two parts, a fluctuation term and an average
term

\begin{align*}
 & -\left(\mathcal{\rho}_{n}\mathcal{S}_{n}\right)^{-1}\sum_{x}\intop_{0}^{t}\xi_{s-}^{n}(x)\phi_{s}(x)(d\hat{P}_{s}(x)-P(\sum_{w\sim x}\xi_{s-}^{n}(w)(\mathcal{N}_{n}\mathcal{S}_{n})^{-1})ds)\\
 & -\left(\mathcal{\rho}_{n}\mathcal{S}_{n}\right)^{-1}\sum_{x}\intop_{0}^{t}\xi_{s-}^{n}(x)\phi_{s}(x)P(\sum_{w\sim x}\xi_{s-}^{n}(w)(\mathcal{N}_{n}\mathcal{S}_{n})^{-1})ds\\
 & =E_{t}^{(2)}(\phi)-\intop_{0}^{t}(v_{s},P(u_{s})\phi_{s})ds,
\end{align*}

where
\begin{align*}
E_{t}^{(2)}(\phi) & :=-\left(\mathcal{\rho}_{n}\mathcal{S}_{n}\right)^{-1}\sum_{x}\intop_{0}^{t}\xi_{s-}^{n}(x)\phi_{s}(x)(d\hat{P}_{s}(x)-P(\sum_{w\sim x}\xi_{s-}^{n}(w)(\mathcal{N}_{n}\mathcal{S}_{n})^{-1})ds)
\end{align*}

is a martingale with brackets process given by

\begin{align*}
d\left\langle E_{n}^{(2)}(\phi)\right\rangle _{t} & =\left(\mathcal{\rho}_{n}\mathcal{S}_{n}\right)^{-2}\sum_{x}\xi_{t}^{n}(x)\phi_{t}^{2}(x)P(\sum_{w\sim x}\xi_{s-}^{n}(w)(\mathcal{N}_{n}\mathcal{S}_{n})^{-1})dt\\
 & =\left(\mathcal{\rho}_{n}\mathcal{S}_{n}\right)^{-1}(v_{t},\phi_{t}^{2}P(u_{t}))dt\\
 & \leqslant\left(\mathcal{\rho}_{n}\mathcal{S}_{n}\right)^{-1}\mathcal{L}_{n}^{c}(1,\phi_{t}^{2}P)dt\\
 & \leqslant\left(\mathcal{\rho}_{n}\mathcal{S}_{n}\right)^{-1}\mathcal{L}_{n}^{c}\left\Vert \phi_{t}^{2}P\right\Vert _{2\lambda}(e_{-2\lambda},1)dt.
\end{align*}

Alternatively, if $-1\leqslant u_{t}\leqslant1$,we bound it by
\begin{align*}
d\left\langle E_{n}^{(2)}(\phi)\right\rangle _{t} & \leqslant\left(\mathcal{\rho}_{n}\mathcal{S}_{n}\right)^{-1}\left|P\right|\mathcal{L}_{n}^{c}(v_{t},\phi_{t}^{2})dt\\
 & \leqslant\left(\mathcal{\rho}_{n}\mathcal{S}_{n}\right)^{-1}\left|P\right|\mathcal{L}_{n}^{c}(1,\phi_{t}^{2})dt\\
 & \leqslant\left(\mathcal{\rho}_{n}\mathcal{S}_{n}\right)^{-1}\left|P\right|\mathcal{L}_{n}^{c}\left\Vert \phi_{t}\right\Vert _{\lambda}^{2}(e_{-2\lambda},1)dt.
\end{align*}
We break the term into two parts, a fluctuation term and an average
term

\begin{align*}
 & -\left(\mathcal{\rho}_{n}\mathcal{S}_{n}\right)^{-1}\sum_{x}\sum_{y\sim x}\sum_{z\sim x}\intop_{0}^{t}\xi_{s-}^{n}(x)\phi_{s}(x)(d\tilde{P}_{s}(x,y,z)-(\mathcal{L}_{n}^{c}/\mathcal{N}_{n}^{2})ds)\\
 & -\left(\mathcal{\rho}_{n}\mathcal{S}_{n}\right)^{-1}\sum_{x}\sum_{y\sim x}\sum_{z\sim x}\intop_{0}^{t}\xi_{s-}^{n}(x)\phi_{s}(x)(\mathcal{L}_{n}^{c}/\mathcal{N}_{n}^{2})ds\\
 & =E_{t}^{(3)}(\phi)-\mathcal{L}_{n}^{c}\intop_{0}^{t}(v_{s},\phi_{s})ds,
\end{align*}

where
\begin{align*}
E_{t}^{(3)}(\phi) & :=-\left(\mathcal{\rho}_{n}\mathcal{S}_{n}\right)^{-1}\sum_{x}\sum_{y\sim x}\sum_{z\sim x}\intop_{0}^{t}\xi_{s-}^{n}(x)\phi_{s}(x)(d\tilde{P}_{s}(x,y,z)-(\mathcal{L}_{n}^{c}/\mathcal{N}_{n}^{2})ds)
\end{align*}

is a martingale with brackets process given by

\begin{align*}
d\left\langle E_{n}^{(3)}(\phi)\right\rangle _{t} & =\left(\mathcal{\rho}_{n}\mathcal{S}_{n}\right)^{-2}\sum_{x}\sum_{y\sim x}\sum_{z\sim x}\xi_{t}^{n}(x)\phi_{t}^{2}(x)(\mathcal{L}_{n}^{c}/\mathcal{N}_{n}^{2})dt\\
 & =\left(\mathcal{\rho}_{n}\mathcal{S}_{n}\right)^{-2}\mathcal{L}_{n}^{c}\sum_{x}\xi_{t}^{n}(x)\phi_{t}^{2}(x)dt\\
 & =\left(\mathcal{\rho}_{n}\mathcal{S}_{n}\right)^{-1}\mathcal{L}_{n}^{c}(v_{t},\phi_{t}^{2})dt\\
 & \leqslant\left(\mathcal{\rho}_{n}\mathcal{S}_{n}\right)^{-1}\mathcal{L}_{n}^{c}(1,\phi_{t}^{2})dt\\
 & \leqslant\left(\mathcal{\rho}_{n}\mathcal{S}_{n}\right)^{-1}\mathcal{L}_{n}^{c}\left\Vert \phi_{t}\right\Vert _{\lambda}^{2}(e_{-2\lambda},1)dt.
\end{align*}
We break the term into two parts, a fluctuation term and an average
term

\begin{align*}
 & \left(\mathcal{\rho}_{n}\mathcal{S}_{n}\right)^{-1}\sum_{x}\sum_{y\sim x}\sum_{z\sim x}\intop_{0}^{t}(\sum_{p(x,y,z)}\xi_{s-}^{n}(x)\xi_{s-}^{n}(y))\phi_{s}(x)(d\tilde{P}_{s}(x,y,z)-(\mathcal{L}_{n}^{c}/\mathcal{N}_{n}^{2})ds)\\
 & +\left(\mathcal{\rho}_{n}\mathcal{S}_{n}\right)^{-1}\sum_{x}\sum_{y\sim x}\sum_{z\sim x}\intop_{0}^{t}(\sum_{p(x,y,z)}\xi_{s-}^{n}(x)\xi_{s-}^{n}(y))\phi_{s}(x)(\mathcal{L}_{n}^{c}/\mathcal{N}_{n}^{2})ds\\
 & =E_{t}^{(4)}(\phi)+2\mathcal{L}_{n}^{c}\intop_{0}^{t}(v_{s},\mathcal{S}_{n}u_{s}\phi_{s})ds+\mathcal{L}_{n}^{c}\intop_{0}^{t}(1,\mathcal{S}_{n}^{2}u_{s}^{2}\phi_{s})ds,
\end{align*}

where
\begin{align*}
E_{t}^{(4)}(\phi) & :=-\left(\mathcal{\rho}_{n}\mathcal{S}_{n}\right)^{-1}\sum_{x}\sum_{y\sim x}\sum_{z\sim x}\intop_{0}^{t}(\sum_{p(x,y,z)}\xi_{s-}^{n}(x)\xi_{s-}^{n}(y))\phi_{s}(x)(d\tilde{P}_{s}(x,y,z)-(\mathcal{L}_{n}^{c}/\mathcal{N}_{n}^{2})ds)
\end{align*}

is a martingale with brackets process given by

\begin{align*}
d\left\langle E_{n}^{(4)}(\phi)\right\rangle _{t} & =\left(\mathcal{\rho}_{n}\mathcal{S}_{n}\right)^{-2}\sum_{x}\sum_{y\sim x}\sum_{z\sim x}(\sum_{p(x,y,z)}\xi_{s-}^{n}(x)\xi_{s-}^{n}(y))^{2}\phi_{s}^{2}(x)(\mathcal{L}_{n}^{c}/\mathcal{N}_{n}^{2})dt\\
 & =\left(\mathcal{\rho}_{n}\mathcal{S}_{n}\right)^{-2}\sum_{x}\sum_{y\sim x}\sum_{z\sim x}[(\sum_{p(x,y,z)}\xi_{s-}^{n}(x)\xi_{s-}^{n}(y))+6\xi_{s-}^{n}(x)\xi_{s-}^{n}(y)\xi_{s-}^{n}(z)]\phi_{s}^{2}(x)(\mathcal{L}_{n}^{c}/\mathcal{N}_{n}^{2})dt\\
 & =\left(\mathcal{\rho}_{n}\mathcal{S}_{n}\right)^{-1}\mathcal{L}_{n}^{c}[2(v_{t},\phi_{t}^{2}u_{t}\mathcal{S}_{n})+(1,u_{t}^{2}\mathcal{S}_{n}^{2}\phi_{t}^{2})+6(v_{t},\phi_{t}^{2}u_{t}^{2}\mathcal{S}_{n}^{2})]dt\\
 & \leqslant9\left(\mathcal{\rho}_{n}\mathcal{S}_{n}\right)^{-1}\mathcal{L}_{n}^{c}(1,\phi_{t}^{2})dt\\
 & \leqslant9\left(\mathcal{\rho}_{n}\mathcal{S}_{n}\right)^{-1}\mathcal{L}_{n}^{c}\left\Vert \phi_{t}\right\Vert _{\lambda}^{2}(e_{-2\lambda},1)dt.
\end{align*}

We break the term into two parts, a fluctuation term and an average
term

\begin{align*}
 & -2\left(\mathcal{\rho}_{n}\mathcal{S}_{n}\right)^{-1}\sum_{x}\sum_{y\sim x}\sum_{z\sim x}\intop_{0}^{t}\xi_{s-}^{n}(x)\xi_{s-}^{n}(y)\xi_{s-}^{n}(z)\phi_{s}(x)(d\tilde{P}_{s}(x,y,z)-(\mathcal{L}_{n}^{c}/\mathcal{N}_{n}^{2})ds)\\
 & -2\left(\mathcal{\rho}_{n}\mathcal{S}_{n}\right)^{-1}\sum_{x}\sum_{y\sim x}\sum_{z\sim x}\intop_{0}^{t}\xi_{s-}^{n}(x)\xi_{s-}^{n}(y)\xi_{s-}^{n}(z)\phi_{s}(x)(\mathcal{L}_{n}^{c}/\mathcal{N}_{n}^{2})ds\\
 & =E_{t}^{(5)}(\phi)-2\mathcal{L}_{n}^{c}\intop_{0}^{t}(v_{s},\mathcal{S}_{n}^{2}u_{s}^{2}\phi_{s})ds,
\end{align*}

where
\begin{align*}
E_{t}^{(5)}(\phi) & :=-2\left(\mathcal{\rho}_{n}\mathcal{S}_{n}\right)^{-1}\sum_{x}\sum_{y\sim x}\sum_{z\sim x}\intop_{0}^{t}\xi_{s-}^{n}(x)\xi_{s-}^{n}(y)\xi_{s-}^{n}(z)\phi_{s}(x)(d\tilde{P}_{s}(x,y,z)-(\mathcal{L}_{n}^{c}/\mathcal{N}_{n}^{2})ds)
\end{align*}

is a martingale with brackets process given by

\begin{align*}
d\left\langle E_{n}^{(5)}(\phi)\right\rangle _{t} & =4\left(\mathcal{\rho}_{n}\mathcal{S}_{n}\right)^{-2}\sum_{x}\sum_{y\sim x}\sum_{z\sim x}\xi_{s-}^{n}(x)\xi_{s-}^{n}(y)\xi_{s-}^{n}(z)\phi_{t}^{2}(x)(\mathcal{L}_{n}^{c}/\mathcal{N}_{n}^{2})dt\\
 & =4\left(\mathcal{\rho}_{n}\mathcal{S}_{n}\right)^{-1}\mathcal{L}_{n}^{c}(v_{t},u_{t}^{2}\mathcal{S}_{n}^{2}\phi_{t}^{2})dt\\
 & \leqslant4\left(\mathcal{\rho}_{n}\mathcal{S}_{n}\right)^{-1}\mathcal{L}_{n}^{c}(1,\phi_{t}^{2})dt\\
 & \leqslant4\left(\mathcal{\rho}_{n}\mathcal{S}_{n}\right)^{-1}\mathcal{L}_{n}^{c}\left\Vert \phi_{t}\right\Vert _{\lambda}^{2}(e_{-2\lambda},1)dt.
\end{align*}

\subsection{White-Noise Term}

We break the term into two parts, a fluctuation term and an average
term

\begin{align*}
 & \left(\mathcal{\rho}_{n}\mathcal{S}_{n}\right)^{-1}\sum_{x}\sum_{y\sim x}\intop_{0}^{t}\xi_{s-}^{n}(x)\phi_{s}(x)(dP_{s}(y,x)-dP_{s}(x,y))\\
 & =\left(\mathcal{\rho}_{n}\mathcal{S}_{n}\right)^{-1}\sum_{x}\sum_{y\sim x}\intop_{0}^{t}\xi_{s-}^{n}(x)\phi_{s}(x)((dP_{s}(y,x)-(\mathcal{H}_{n}/\mathcal{N}_{n})ds)-(dP_{s}(x,y)-(\mathcal{H}_{n}/\mathcal{N}_{n})ds))\\
 & =Z_{t}(\phi),
\end{align*}

where $Z_{t}(\phi)$ is a martingale with brackets process given by

\begin{align}
d\left\langle Z(\phi)\right\rangle _{t} & =\left(\mathcal{\rho}_{n}\mathcal{S}_{n}\right)^{-2}\sum_{x}\sum_{y\sim x}\sum_{x^{\prime}}\sum_{y^{\prime}\sim x^{\prime}}\xi_{t}^{n}(x)\phi_{t}(x)\xi_{t}^{n}(x^{\prime})\phi_{t}(x^{\prime})2I(x=x^{\prime},y=y^{\prime})(\mathcal{H}_{n}/\mathcal{N}_{n})dt\nonumber \\
 & -\left(\mathcal{\rho}_{n}\mathcal{S}_{n}\right)^{-2}\sum_{x}\sum_{y\sim x}\sum_{x^{\prime}}\sum_{y^{\prime}\sim x^{\prime}}\xi_{t}^{n}(x)\phi_{t}(x)\xi_{t}^{n}(x^{\prime})\phi_{t}(x^{\prime})2I(x=y^{\prime},y=x^{\prime})(\mathcal{H}_{n}/\mathcal{N}_{n})dt\nonumber \\
 & =2\left(\mathcal{\rho}_{n}\mathcal{S}_{n}\right)^{-2}\sum_{x}\sum_{y\sim x}(\xi_{t}^{n}(x)\phi_{t}(x)^{2}-\xi_{t}^{n}(x)\phi_{t}(x)\xi_{t}^{n}(y)\phi_{t}(y))(\mathcal{H}_{n}/\mathcal{N}_{n})dt\nonumber \\
 & =2\left(\mathcal{\rho}_{n}\mathcal{S}_{n}\right)^{-1}\mathcal{H}_{n}\left[(v_{t},\phi_{t}^{2})-(v_{t},\mathcal{S}_{n}\phi_{t}A(\xi_{t}^{n}\phi_{t}))\right]dt\nonumber \\
 & \leqslant2\left(\mathcal{\rho}_{n}\mathcal{S}_{n}\right)^{-1}\mathcal{H}_{n}(v_{t},\phi_{t}^{2})dt\nonumber \\
 & \leqslant2\left(\mathcal{\rho}_{n}\mathcal{S}_{n}\right)^{-1}\mathcal{H}_{n}\left\Vert \phi_{t}\right\Vert _{\lambda}^{2}(e_{-2\lambda},1)dt.\label{eq:2}
\end{align}

Collecting all terms above, we get the final semi-martingale decomposition:

\begin{align*}
(v_{t},\phi_{t}) & =(v_{0},\phi_{0})+\intop_{0}^{t}\left(v_{s},\partial_{s}\phi_{s}+\Delta_{n}(\phi_{s})-\mathcal{L}_{n}^{c}\phi_{s}\right)ds\\
 & +\intop_{0}^{t}\left(v_{s},-P(u_{s})\phi_{s}+2\mathcal{L}_{n}^{c}\mathcal{S}_{n}u_{s}\phi_{s}-2\mathcal{L}_{n}^{c}\mathcal{S}_{n}^{2}u_{s}^{2}\phi_{s}\right)ds+\intop_{0}^{t}(1,\mathcal{L}_{n}^{c}\mathcal{S}_{n}^{2}u_{s}^{2}\phi_{s})ds\\
 & +\sum_{i=1}^{5}E_{t}^{(i)}(\phi)+Z_{t}(\phi).
\end{align*}

\section{Green's Function Representation}

We need to choose a special test function to the second term in the
above final semi-martingale decomposition. For each $z\in\mathbb{Z}/\mathcal{\rho}_{n}$,
we define:
\begin{enumerate}
\item a function $\psi_{t}^{z}$ to be the unique solution of 
\begin{align*}
\partial_{t}\psi_{t}^{z} & =\Delta_{n}\psi_{t}^{z}\\
\psi_{0}^{z}(x) & =\mathcal{\rho}_{n}\mathcal{N}_{n}^{-1}I(x\sim z).
\end{align*}
\item a random walk $X_{t}$ with generator $\Delta_{n}$ jumping at rate
$\mathcal{H}_{n}$ with symmetric steps of variance $\mathcal{D}_{n}^{2}/3.$
\item $\bar{\psi}_{0}^{z}(x):=\mathcal{\rho}_{n}P(X_{t}=x|X_{0}=z).$ 
\end{enumerate}
Maybe you note there is a relation between $\psi_{t}^{z}$ and $\bar{\psi}_{t}^{z}$,
i.e. 
\[
\psi_{t}^{z}(x)=\underset{w\sim z}{\sum}\mathcal{\rho}_{n}\mathcal{N}_{n}^{-1}P(X_{t}=x|X_{0}=w)=\underset{w}{\sum}\mathcal{\rho}_{n}\mathcal{N}_{n}^{-1}I(w\sim z)P(X_{t}=w|X_{0}=x)=(\psi_{0}^{z},\bar{\psi}_{t}^{x}),
\]

and a property that 
\[
\psi_{t}^{z}(x)=\psi_{t}^{x}(z)\quad and\quad\bar{\psi}_{t}^{z}(x)=\bar{\psi}_{t}^{x}(z).
\]

Then linearly interpolate $\psi_{t}^{z}$ and $\bar{\psi}_{t}^{z}$,
you will get a convergence of them to $p(\mathcal{D}_{n}^{2}\mathcal{H}_{n}t/3,x-z)$.

\subsection{Property of $\psi_{t}^{z}$ and $\bar{\psi}_{t}^{z}$}

Under different coefficients, we have similar estimations associated
with $\psi_{t}^{z}$ and $\bar{\psi}_{t}^{z}$ to \cite[Lemma 3]{mueller1995stochastic}.
In our case we can use the conclusion.

Set $\phi_{s}=e^{-\mathcal{L}_{n}^{c}(t-s)}\psi_{t-s}^{x}$ for $s\leqslant t$
and substitute it into the final semi-martingale decomposition equation
then the second term vanishes and $(v_{t},\phi_{t})=u_{t}(x)$, i.e.
$u_{n}(t,x)$. Hence it turns out that:

\begin{align*}
u_{t}(x) & =(v_{0},e^{-\mathcal{L}_{n}^{c}t}\psi_{t}^{x})+\intop_{0}^{t}\left(v_{s},-P(u_{s})e^{-\mathcal{L}_{n}^{c}(t-s)}\psi_{t-s}^{x}\right)ds\\
 & +2\intop_{0}^{t}\left(v_{s},\mathcal{L}_{n}^{c}\mathcal{S}_{n}u_{s}e^{-\mathcal{L}_{n}^{c}(t-s)}\psi_{t-s}^{x}-\mathcal{L}_{n}^{c}\mathcal{S}_{n}^{2}u_{s}^{2}e^{-\mathcal{L}_{n}^{c}(t-s)}\psi_{t-s}^{x}\right)ds\\
 & +\intop_{0}^{t}(1,\mathcal{L}_{n}^{c}\mathcal{S}_{n}^{2}u_{s}^{2}e^{-\mathcal{L}_{n}^{c}(t-s)}\psi_{t-s}^{x})ds+\sum_{i=1}^{5}E_{t}^{(i)}(e^{-\mathcal{L}_{n}^{c}(t-s)}\psi_{t-s}^{x})+Z_{t}(e^{-\mathcal{L}_{n}^{c}(t-s)}\psi_{t-s}^{x}).
\end{align*}

Because we will only use this representation of $u_{n}(t,x)$ to estimate
moment required for the proof of tightness ($e^{-\mathcal{L}_{n}^{c}(t-s)}$
is bounded and positive), we reduce it to 
\begin{align}
u_{t}(x) & =(v_{0},\psi_{t}^{x})+\intop_{0}^{t}\left(v_{s},-P(u_{s})\psi_{t-s}^{x}+2\mathcal{L}_{n}^{c}\mathcal{S}_{n}u_{s}\psi_{t-s}^{x}-2\mathcal{L}_{n}^{c}\mathcal{S}_{n}^{2}u_{s}^{2}\psi_{t-s}^{x}\right)ds\label{eq:1}\\
 & +\intop_{0}^{t}(1,\mathcal{L}_{n}^{c}\mathcal{S}_{n}^{2}u_{s}^{2}\psi_{t-s}^{x})ds+\sum_{i=1}^{5}E_{t}^{(i)}(\psi_{t-\cdot}^{x})+Z_{t}(\psi_{t-\cdot}^{x}),\nonumber 
\end{align}

for convenience without loss of generality.

If $\mathcal{L}_{n}^{c},\mathcal{S}_{n}$ are constant, for simplicity,
we rewrite the above formula:
\[
u_{t}(x)=(v_{0},\psi_{t}^{x})+\intop_{0}^{t}\left(v_{s},\mathcal{P}(u_{s})\psi_{t-s}^{x}\right)ds+\intop_{0}^{t}(1,\mathcal{L}_{n}^{c}\mathcal{S}_{n}^{2}u_{s}^{2}\psi_{t-s}^{x})ds+\sum_{i=1}^{5}E_{t}^{(i)}(\psi_{t-\cdot}^{x})+Z_{t}(\psi_{t-\cdot}^{x}),
\]

where 
\[
\hat{P}(u_{s})=-P(u_{s})+2\mathcal{L}_{n}^{c}\mathcal{S}_{n}u_{s}-2\mathcal{L}_{n}^{c}\mathcal{S}_{n}^{2}u_{s}^{2}.
\]

\begin{lem}
If $f:\mathbb{Z}/\rho_{n}\longrightarrow[0,\infty)$ with $(f,f)<\infty,\lambda\in\mathbb{R},m\in\mathbb{N^{+}},$
then:
\end{lem}

\begin{enumerate}
\item $(v_{t},\psi_{s}^{z})=(u_{t},\bar{\psi}_{s}^{z}),$
\item $\left|(v_{t},f)-(u_{t},f)\right|\leqslant\left\Vert D(f,\mathcal{D}_{n})\right\Vert _{\lambda}(v_{t},e_{-\lambda}),$
\item $(v_{t},u_{t}^{m}e_{-\lambda})\leqslant C(m)(u_{t}^{m+1},e_{-\lambda})$
(if $u_{t}(x)$ is bounded).
\end{enumerate}
Proof:

\begin{align*}
(v_{t},\psi_{s}^{z}) & =(\mathcal{\rho}_{n}\mathcal{S}_{n})^{-1}\underset{x}{\sum}\xi_{t}(x)\underset{w\sim z}{\sum}\mathcal{\rho}_{n}\mathcal{N}_{n}^{-1}P(X_{s}=x|X_{0}=w)\\
 & =(\mathcal{\rho}_{n}\mathcal{S}_{n})^{-1}\underset{x}{\sum}\xi_{t}(x)\underset{w}{\sum}\mathcal{\rho}_{n}\mathcal{N}_{n}^{-1}P(X_{s}=z|X_{0}=w)I(w\sim x)\\
 & =\mathcal{\rho}_{n}^{-1}\underset{x}{\sum}\xi_{t}(x)\underset{w}{\sum}\mathcal{\rho}_{n}(\mathcal{S}_{n}\mathcal{N}_{n})^{-1}P(X_{s}=w|X_{0}=z)I(w\sim x)\\
 & =(u_{t},\bar{\psi}_{s}^{z}).
\end{align*}

\begin{align*}
(u_{t},f) & =\mathcal{\rho}_{n}^{-1}(\mathcal{S}_{n}\mathcal{N}_{n})^{-1}\underset{x}{\sum}\sum_{y\sim x}\xi_{t}(y)f(x)\\
 & =\mathcal{\rho}_{n}^{-1}(\mathcal{S}_{n}\mathcal{N}_{n})^{-1}\underset{x}{\sum}\sum_{y\sim x}\xi_{t}(y)(f(x)-f(y))+(v_{t},f),
\end{align*}
so
\begin{align*}
\left|(v_{t},f)-(u_{t},f)\right| & =\mathcal{\rho}_{n}^{-1}(\mathcal{S}_{n}\mathcal{N}_{n})^{-1}\underset{x}{\sum}\sum_{y\sim x}\xi_{t}(y)\left|f(x)-f(y)\right|\\
 & \leqslant(v_{t},D(f,\mathcal{D}_{n}))\\
 & \leqslant\left\Vert D(f,\mathcal{D}_{n})\right\Vert _{\lambda}(v_{t},e_{-\lambda}).
\end{align*}

\begin{align*}
(u_{t}^{2},f) & =\mathcal{\rho}_{n}^{-1}\underset{x}{\sum}f(x)(\sum_{y\sim x}\xi_{t}(y)(\mathcal{S}_{n}\mathcal{N}_{n})^{-1})^{2}\\
 & \geqslant\mathcal{\rho}_{n}^{-1}\underset{x}{\sum}\sum_{y\sim x}(f(y)-D(f,\mathcal{D}_{n})(x))\xi_{t}(y)(\mathcal{S}_{n}\mathcal{N}_{n})^{-1}(\sum_{w\sim x}\xi_{t}(w)(\mathcal{S}_{n}\mathcal{N}_{n})^{-1})\\
 & \geqslant\mathcal{\rho}_{n}^{-1}\underset{x}{\sum}\sum_{y\sim x}f(y)\xi_{t}(y)(\mathcal{S}_{n}\mathcal{N}_{n})^{-1}(\sum_{w\sim x}\xi_{t}(w)(\mathcal{S}_{n}\mathcal{N}_{n})^{-1})-(u_{t}^{2},D(f,\mathcal{D}_{n}))\\
 & \geqslant\mathcal{\rho}_{n}^{-1}\underset{y}{\sum}f(y)\xi_{t}(y)(\mathcal{S}_{n}\mathcal{N}_{n})^{-1}(\sum_{w\sim y}\xi_{t}(w)(\mathcal{S}_{n}\mathcal{N}_{n})^{-1})\sum_{x}I(x\sim w,x\sim y)-(u_{t}^{2},D(f,\mathcal{D}_{n}))\\
 & \geqslant C(\mathcal{S}_{n}\mathcal{\rho}_{n})^{-1}\underset{y}{\sum}f(y)\xi_{t}(y)(\sum_{w\sim y}\xi_{t}(w)(\mathcal{S}_{n}\mathcal{N}_{n})^{-1})-(u_{t}^{2},D(f,\mathcal{D}_{n}))\\
 & \geqslant C(v_{t},u_{t}f)-(u_{t}^{2},D(f,\mathcal{D}_{n})),
\end{align*}

i.e. 
\[
(v_{t},u_{t}f)\leqslant C(u_{t}^{2},f)+C(u_{t}^{2},D(f,\mathcal{D}_{n}))
\]

Then substitute $u_{t}^{m-1}e_{-\lambda}$ for $f$, we obtain:
\begin{align*}
(v_{t},u_{t}^{m}e_{-\lambda}) & \leqslant C(u_{t}^{m+1},e_{-\lambda})+C(u_{t}^{2},D(u_{t}^{m-1}e_{-\lambda},\mathcal{D}_{n}))\\
 & \leqslant C(m)(u_{t}^{m+1},e_{-\lambda}).
\end{align*}
From now on, the incoming method is valid for the cases where $\mathcal{S}_{n}=1,0\leqslant u_{n}(0,x)\leqslant1$,
since the method branches according to $\mathcal{S}_{n}$ with different
value, however, we still use notation $\mathcal{S}_{n}$ to identify
its trace.
\begin{lem}
For $T\geqslant0,\;p\geqslant2,\lambda>0,$ we have:
\end{lem}

\[
E(\left|E_{t}^{(i)}(\psi_{t-\cdot}^{z})\right|^{p})\leqslant C(\lambda,p,T)n^{-p/16}e_{\lambda p}(z)\quad for\;all\;t\leqslant T.
\]
Proof:

For $i\geqslant2$, we can prove it in an identical way. The greatest
jump of these martingale is bounded by $\left(\mathcal{\rho}_{n}\mathcal{S}_{n}\right)^{-1}\underset{s\geqslant0}{\sup}\left\Vert \psi_{s}^{z}\right\Vert _{\infty}\leqslant C\left(\mathcal{\rho}_{n}\mathcal{S}_{n}\right)^{-1}n^{1/2}$.
So by implementing Burkholder's inequality, we get:

\begin{align*}
 & E(\left|E_{t}^{(i)}(\psi_{t-\cdot}^{z})\right|^{p})\\
 & \leqslant C(\lambda,p,T)\left(\mathcal{\rho}_{n}^{-1}\mathcal{S}_{n}^{-1}\right)^{p/2}((\mathcal{L}_{n}^{c}\intop_{0}^{t}\left\Vert \psi_{t-s}^{z}\right\Vert _{\lambda}^{2}ds)^{p/2}+\left(\mathcal{\rho}_{n}^{-1}\mathcal{S}_{n}^{-1}n\right)^{p/2})\\
 & \leqslant C(\lambda,p,T)\left(\mathcal{\rho}_{n}^{-1}\mathcal{S}_{n}^{-1}\right)^{p/2}((\mathcal{L}_{n}^{c}n^{1/4})^{p/2}+\left(\mathcal{\rho}_{n}^{-1}\mathcal{S}_{n}^{-1}n\right)^{p/2})e_{\lambda p}(z).
\end{align*}

For $i=1$, also use Burkholder's inequality, we have:
\begin{align*}
 & E(\left|E_{t}^{(1)}(\psi_{t-\cdot}^{z})\right|^{p})\\
 & \leqslant C(p)\left(\mathcal{\rho}_{n}^{-1}\mathcal{S}_{n}^{-1}\mathcal{H}_{n}\right)^{p/2}E((\intop_{(t-n^{-3/8})_{+}}^{t}\left\Vert \psi_{t-s}^{z}\right\Vert _{0}\left[(\psi_{t-s}^{z},u_{t}\mathcal{S}_{n})+(v_{t},\psi_{t-s}^{z})\right]ds)^{p/2})\\
 & +C(p)\left(\mathcal{\rho}_{n}^{-1}\mathcal{S}_{n}^{-1}\mathcal{H}_{n}\right)^{p/2}E((\intop_{0}^{(t-n^{-3/8})_{+}}\left\Vert D(\psi_{t-s}^{z},\mathcal{D}_{n})\right\Vert _{\lambda}^{2}(e_{-2\lambda},v_{s})ds)^{p/2})\\
 & +C(p)\left(\mathcal{\rho}_{n}\mathcal{S}_{n}\right)^{-p}n^{p/2}
\end{align*}

Let us have a look at the first term:

\begin{align*}
 & E((\intop_{(t-n^{-3/8})_{+}}^{t}\left\Vert \psi_{t-s}^{z}\right\Vert _{0}\left[(\psi_{t-s}^{z},u_{t}\mathcal{S}_{n})+(v_{t},\psi_{t-s}^{z})\right]ds)^{p/2})\\
 & \leqslant C(T)E((\intop_{(t-n^{-3/8})_{+}}^{t}(t-s)^{-2/3}(u_{t},\psi_{t-s}^{z}+\bar{\psi}_{t-s}^{z})ds)^{p/2})\\
 & \leqslant C(T)(\intop_{(t-n^{-3/8})_{+}}^{t}(t-s)^{-2/3}(1,\psi_{t-s}^{z}+\bar{\psi}_{t-s}^{z})ds)^{p/2}\\
 & \leqslant C(\lambda,p,T)n^{-p/16}e_{\lambda p}(z),
\end{align*}

then the second term:
\begin{align*}
 & E((\intop_{0}^{(t-n^{-3/8})_{+}}\left\Vert D(\psi_{t-s}^{z},\mathcal{D}_{n})\right\Vert _{\lambda}^{2}(e_{-2\lambda},v_{s})ds)^{p/2})\\
 & \leqslant C(\lambda,p,T)(\intop_{0}^{(t-n^{-3/8})_{+}}\left\Vert D(\psi_{t-s}^{z},\mathcal{D}_{n})\right\Vert _{\lambda}^{2}ds)^{p/2}\\
 & \leqslant C(\lambda,p,T)e_{\lambda p}(z)n^{-p/4}(\intop_{0}^{(t-n^{-3/8})_{+}}(t-s)^{-2}ds)^{p/2}\\
 & \leqslant C(\lambda,p,T)n^{-p/16}e_{\lambda p}(z).
\end{align*}

Finally, we get:
\[
E(\left|E_{t}^{(1)}(\psi_{t-\cdot}^{z})\right|^{p})\leqslant C(\lambda,p,T)e_{\lambda p}(z)\left(\mathcal{\rho}_{n}^{-1}\mathcal{S}_{n}^{-1}\right)^{p/2}(\left(\mathcal{H}_{n}n^{-1/8}\right)^{p/2}+\left(\mathcal{\rho}_{n}^{-1}\mathcal{S}_{n}^{-1}n\right)^{p/2}).
\]

According to our assumption at the beginning, our version of the above
bound is:
\[
E(\left|E_{t}^{(i)}(\psi_{t-\cdot}^{z})\right|^{p})\leqslant C(\lambda,p,T)n^{-p/16}e_{\lambda p}(z).
\]

\section{Tightness}

Our objective is to prove the tightness of $u_{n}(t,x)$ by estimating
moment of their discrepancy at different times and locations.
\begin{lem}
For $0\leqslant s\leqslant t\leqslant T,\;y,z\in\mathbb{Z}/\mathcal{\rho}_{n},\;\left|t-s\right|,\left|y-z\right|\leqslant1,\;\lambda>0,\;p\geqslant2,\;\mathcal{L}_{n}^{c},\mathcal{S}_{n}\;are\;constant$
\end{lem}

\[
E(\left|\hat{u}_{t}(z)-\hat{u}_{s}(y)\right|^{p})\leqslant C(\lambda,p,T)e_{\lambda p}(z)(\left|t-s\right|^{p/24}+\left|z-y\right|^{p/24}+n^{-p/12}).
\]
Proof:
\[
\hat{u}_{t}(x)=\intop_{0}^{t}\left(v_{s},\hat{P}(u_{s})\psi_{t-s}^{x}\right)ds+\intop_{0}^{t}(1,\mathcal{L}_{n}^{c}\mathcal{S}_{n}^{2}u_{s}^{2}\psi_{t-s}^{x})ds+\sum_{i=1}^{5}E_{t}^{(i)}(\psi_{t-\cdot}^{x})+Z_{t}(\psi_{t-\cdot}^{x}).
\]

\[
E(\left|\hat{u}_{t}(z)-\hat{u}_{s}(y)\right|^{p})\leqslant C(p)E(\left|\hat{u}_{t}(z)-\hat{u}_{t}(y)\right|^{p})+C(p)E(\left|\hat{u}_{t}(z)-\hat{u}_{s}(z)\right|^{p})
\]

Set
\begin{align*}
\delta & =(\left|z-y\right|^{1/4}\lor n^{-1/2})\land t\\
\bar{\delta} & =(\left|t-s\right|^{1/4}\lor n^{-1/2})\land s
\end{align*}

First, we pay attention to the first term:
\begin{align*}
 & E(\left|\hat{u}_{t}(z)-\hat{u}_{t}(y)\right|^{p})\\
 & \leqslant C(\lambda,p,T)n^{-p/16}e_{\lambda p}(z)+C(p)E\left|\intop_{0}^{t}\left(v_{s},\hat{P}(u_{s})(\psi_{t-s}^{z}-\psi_{t-s}^{y})\right)ds\right|^{p}\\
 & +C(p)E\left|\intop_{0}^{t}(1,u_{s}^{2}(\psi_{t-s}^{z}-\psi_{t-s}^{y}))ds\right|^{p}+C(p)E\left|Z_{t}(\psi_{t-\cdot}^{z}-\psi_{t-\cdot}^{y})\right|^{p}\\
 & \leqslant C(\lambda,p,T)n^{-p/16}e_{\lambda p}(z)\\
 & +C(p)(E\left|\intop_{0}^{t-\delta}\left(v_{s},\hat{P}(u_{s})(\psi_{t-s}^{z}-\psi_{t-s}^{y})\right)ds\right|^{p}+E\left|\intop_{t-\delta}^{t}\left(v_{s},\hat{P}(u_{s})(\psi_{t-s}^{z}-\psi_{t-s}^{y})\right)ds\right|^{p})\\
 & +C(p)(E\left|\intop_{0}^{t-\delta}(1,u_{s}^{2}(\psi_{t-s}^{z}-\psi_{t-s}^{y}))ds\right|^{p}+E\left|\intop_{t-\delta}^{t}(1,u_{s}^{2}(\psi_{t-s}^{z}-\psi_{t-s}^{y}))ds\right|^{p})\\
 & +C(p)\left(\mathcal{\rho}_{n}^{-1}\mathcal{S}_{n}^{-1}\mathcal{H}_{n}\right)^{p/2}(E\left|\intop_{0}^{t-\delta}\left(v_{s},(\psi_{t-s}^{z}-\psi_{t-s}^{y})^{2}\right)ds\right|^{p/2}+E\left|\intop_{t-\delta}^{t}\left(v_{s},(\psi_{t-s}^{z}-\psi_{t-s}^{y})^{2}\right)ds\right|^{p/2})
\end{align*}

Let us look at the 1st expectation:
\begin{align*}
 & E\left|\intop_{0}^{t-\delta}\left(v_{s},\hat{P}(u_{s})(\psi_{t-s}^{z}-\psi_{t-s}^{y})\right)ds\right|^{p}\\
 & \leqslant E\left|\intop_{0}^{t-\delta}\left(v_{s},\hat{P}(u_{s})e_{-\lambda}\right)ds\right|^{p}\sup_{s\in(\delta,t]}\left\Vert \psi_{s}^{z}-\psi_{s}^{y}\right\Vert _{\lambda}^{p}\\
 & \leqslant C(\lambda,p,T)e_{\lambda p}(z)(\left|z-y\right|^{p/2}\delta^{-p}+n^{-p/2}\delta^{-3p/4}).
\end{align*}

Let us look at the 2nd expectation:
\begin{align*}
 & E\left|\intop_{t-\delta}^{t}\left(v_{s},\hat{P}(u_{s})(\psi_{t-s}^{z}-\psi_{t-s}^{y})\right)ds\right|^{p}\\
 & \leqslant C(\lambda,p,T)E\left|\intop_{t-\delta}^{t}\left\Vert \psi_{t-s}^{z}\right\Vert _{\lambda}\left(v_{s},\hat{P}(u_{s})e_{-\lambda}\right)ds\right|^{p}\\
 & \leqslant C(\lambda,p,T)E\left|\intop_{t-\delta}^{t}\left\Vert \psi_{t-s}^{z}\right\Vert _{\lambda}\left(u_{s}\hat{P}(u_{s}),e_{-\lambda}\right)ds\right|^{p}\\
 & \leqslant C(\lambda,p,T)e_{\lambda p}(z)\delta^{(1/3)(p-1)}\intop_{t-\delta}^{t}\left|t-s\right|^{-2/3}\left(1,e_{-\lambda}\right)ds\\
 & \leqslant C(\lambda,p,T)\delta^{p/3}e_{\lambda p}(z).
\end{align*}

Let us look at the 3rd expectation:
\begin{align*}
 & E\left|\intop_{0}^{t-\delta}(1,u_{s}^{2}(\psi_{t-s}^{z}-\psi_{t-s}^{y}))ds\right|^{p}\\
 & \leqslant E\left|\intop_{0}^{t-\delta}\left(1,u_{s}^{2}e_{-\lambda}\right)ds\right|^{p}\sup_{s\in(\delta,t]}\left\Vert \psi_{s}^{z}-\psi_{s}^{y}\right\Vert _{\lambda}^{p}\\
 & \leqslant C(\lambda,p,T)e_{\lambda p}(z)(\left|z-y\right|^{p/2}\delta^{-p}+n^{-p/2}\delta^{-3p/4}).
\end{align*}

Let us look at the 4th expectation:
\begin{align*}
 & E\left|\intop_{t-\delta}^{t}(1,u_{s}^{2}(\psi_{t-s}^{z}-\psi_{t-s}^{y}))ds\right|^{p}\\
 & \leqslant C(\lambda,p,T)E\left|\intop_{t-\delta}^{t}\left\Vert \psi_{t-s}^{z}\right\Vert _{\lambda}\left(u_{s}^{2},e_{-\lambda}\right)ds\right|^{p}\\
 & \leqslant C(\lambda,p,T)e_{\lambda p}(z)\delta^{(1/3)(p-1)}\intop_{t-\delta}^{t}\left|t-s\right|^{-2/3}\left(1,e_{-\lambda}\right)ds\\
 & \leqslant C(\lambda,p,T)\delta^{p/3}e_{\lambda p}(z).
\end{align*}

Let us look at the 5th expectation:
\begin{align*}
 & E\left|\intop_{0}^{t-\delta}\left(v_{s},(\psi_{t-s}^{z}-\psi_{t-s}^{y})^{2}\right)ds\right|^{p/2}\\
 & \leqslant E\left|\intop_{0}^{t-\delta}\left(v_{s},e_{-2\lambda}\right)ds\right|^{p/2}\sup_{s\in(\delta,t]}\left\Vert \psi_{s}^{z}-\psi_{s}^{y}\right\Vert _{\lambda}^{p}\\
 & \leqslant C(\lambda,p,T)e_{\lambda p}(z)(\left|z-y\right|^{p/2}\delta^{-p}+n^{-p/2}\delta^{-3p/4}).
\end{align*}

Let us look at the 6th expectation:
\begin{align*}
 & E\left|\intop_{0}^{t-\delta}\left(v_{s},(\psi_{t-s}^{z}-\psi_{t-s}^{y})^{2}\right)ds\right|^{p/2}\\
 & \leqslant E\left|\intop_{t-\delta}^{t}\left\Vert \psi_{t-s}^{z}+\psi_{t-s}^{y}\right\Vert _{0}\left(v_{s},\psi_{t-s}^{z}+\psi_{t-s}^{y}\right)ds\right|^{p/2}\\
 & \leqslant C(T)E\left|\intop_{t-\delta}^{t}\left|t-s\right|^{-2/3}\left(u_{s},\bar{\psi}_{t-s}^{z}+\bar{\psi}_{t-s}^{y}\right)ds\right|^{p/2}\\
 & \leqslant C(\lambda,p,T)\delta^{(1/3)((p/2)-1)}\intop_{t-\delta}^{t}\left|t-s\right|^{-2/3}\left(1,\bar{\psi}_{t-s}^{z}+\bar{\psi}_{t-s}^{y}\right)ds\\
 & \leqslant C(\lambda,p,T)\delta^{p/6}.
\end{align*}

Put all the above conclusions together:
\[
E(\left|\hat{u}_{t}(z)-\hat{u}_{t}(y)\right|^{p})\leqslant C(\lambda,p,T)e_{\lambda p}(z)(\left|z-y\right|^{p/24}+n^{-p/12}).
\]

Similarly, we consider the second term:
\begin{align*}
 & E(\left|\hat{u}_{t}(z)-\hat{u}_{s}(z)\right|^{p}\\
 & \leqslant C(\lambda,p,T)n^{-p/16}e_{\lambda p}(z)\\
 & +C(p)(E\left|\intop_{s}^{t}\left(v_{r},\hat{P}(u_{r})\psi_{t-r}^{z}\right)dr\right|^{p}+E\left|\intop_{0}^{s}\left(v_{r},\hat{P}(u_{r})(\psi_{t-r}^{z}-\psi_{s-r}^{z})\right)dr\right|^{p})\\
 & +C(p)(E\left|\intop_{s}^{t}(1,u_{r}^{2}\psi_{t-r}^{z})dr\right|^{p}+C(p)E\left|\intop_{0}^{s}(1,u_{r}^{2}(\psi_{t-r}^{z}-\psi_{s-r}^{z}))dr\right|^{p})\\
 & +C(p)(E\left|\intop_{s}^{t}\left(v_{r},(\psi_{t-r}^{z})^{2}\right)dr\right|^{p/2}+C(p)E\left|\intop_{0}^{s}\left(v_{r},(\psi_{t-r}^{z}-\psi_{s-r}^{z})^{2}\right)dr\right|^{p/2})\\
 & \leqslant C(\lambda,p,T)n^{-p/16}e_{\lambda p}(z)\\
 & +C(p)(E\left|\intop_{s}^{t}\left(v_{r},\hat{P}(u_{r})\psi_{t-r}^{z}\right)dr\right|^{p}+E\left|\intop_{s}^{t}(1,u_{r}^{2}\psi_{t-r}^{z})dr\right|^{p}+E\left|\intop_{s}^{t}\left(v_{r},(\psi_{t-r}^{z})^{2}\right)dr\right|^{p/2})\\
 & +C(p)(E\left|\intop_{0}^{s-\bar{\delta}}\left(v_{r},\hat{P}(u_{r})(\psi_{t-r}^{z}-\psi_{s-r}^{z})\right)dr\right|^{p}+E\left|\intop_{s-\bar{\delta}}^{s}\left(v_{r},\hat{P}(u_{r})(\psi_{t-r}^{z}-\psi_{s-r}^{z})\right)dr\right|^{p})\\
 & +C(p)(E\left|\intop_{0}^{s-\bar{\delta}}(1,u_{r}^{2}(\psi_{t-r}^{z}-\psi_{s-r}^{z}))dr\right|^{p}+E\left|\intop_{s-\bar{\delta}}^{s}(1,u_{r}^{2}(\psi_{t-r}^{z}-\psi_{s-r}^{z}))dr\right|^{p})\\
 & +C(p)\left(\mathcal{\rho}_{n}^{-1}\mathcal{S}_{n}^{-1}\mathcal{H}_{n}\right)^{p/2}(E\left|\intop_{0}^{s-\bar{\delta}}\left(v_{r},(\psi_{t-r}^{z}-\psi_{s-r}^{z})^{2}\right)dr\right|^{p/2}+E\left|\intop_{s-\bar{\delta}}^{s}\left(v_{r},(\psi_{t-r}^{z}-\psi_{s-r}^{z})^{2}\right)dr\right|^{p/2})
\end{align*}

Let us look at the 1st expectation:
\begin{align*}
 & E\left|\intop_{s}^{t}\left(v_{r},\hat{P}(u_{r})\psi_{t-r}^{z}\right)dr\right|^{p}\\
 & \leqslant E\left|\intop_{s}^{t}\left(v_{r},\hat{P}(u_{r})e_{-\lambda}\right)dr\right|^{p}\sup_{s\leqslant t}\left\Vert \psi_{s}^{z}\right\Vert _{\lambda}^{p}\\
 & \leqslant C(\lambda,p,T)e_{\lambda p}(z)\left|t-s\right|^{p}.
\end{align*}

Let us look at the 2nd expectation:
\begin{align*}
 & E\left|\intop_{s}^{t}(1,u_{r}^{2}\psi_{t-r}^{z})dr\right|^{p}\\
 & \leqslant E\left|\intop_{s}^{t}(1,u_{r}^{2}e_{-\lambda})dr\right|^{p}\sup_{s\leqslant t}\left\Vert \psi_{s}^{z}\right\Vert _{\lambda}^{p}\\
 & \leqslant C(\lambda,p,T)e_{\lambda p}(z)\left|t-s\right|^{p}.
\end{align*}

Let us look at the 3rd expectation:
\begin{align*}
 & E\left|\intop_{s}^{t}\left(v_{r},(\psi_{t-r}^{z})^{2}\right)dr\right|^{p/2}\\
 & \leqslant E\left|\intop_{s}^{t}\left\Vert \psi_{t-r}^{z}\right\Vert _{0}\left(v_{r},(\psi_{t-r}^{z})\right)dr\right|^{p/2}\\
 & \leqslant C(p,T)E\left|\intop_{s}^{t}(t-r)^{-2/3}\left(u_{r},\bar{\psi}_{t-r}^{z}\right)dr\right|^{p/2}\\
 & \leqslant C(p,T)\left|t-s\right|^{p/6}.
\end{align*}

Let us look at the 4th expectation:
\begin{align*}
 & E\left|\intop_{0}^{s-\bar{\delta}}\left(v_{r},\hat{P}(u_{r})(\psi_{t-r}^{z}-\psi_{s-r}^{z})\right)dr\right|^{p}\\
 & \leqslant E\left|\intop_{0}^{s}\left(v_{r},\hat{P}(u_{r})e_{-\lambda}\right)dr\right|^{p}\sup_{r\in[0,s-\bar{\delta})}\left\Vert \psi_{t-r}^{z}-\psi_{s-r}^{z}\right\Vert _{\lambda}^{p}\\
 & \leqslant C(\lambda,p,T)e_{\lambda p}(z)(\left|t-s\right|^{p/2}\bar{\delta}^{-3p/2}+n^{-p/2}\bar{\delta}^{-3p/4}).
\end{align*}

Let us look at the 5th expectation:
\begin{align*}
 & E\left|\intop_{s-\bar{\delta}}^{s}\left(v_{r},\hat{P}(u_{r})(\psi_{t-r}^{z}-\psi_{s-r}^{z})\right)dr\right|^{p}\\
 & \leqslant E\left|\intop_{s-\bar{\delta}}^{s}\left\Vert \psi_{s-r}^{z}\right\Vert _{\lambda}\left(v_{r},\hat{P}(u_{r})e_{-\lambda}\right)dr\right|^{p}\\
 & \leqslant C(\lambda,p,T)e_{\lambda p}(z)E\left|\intop_{s-\bar{\delta}}^{s}\left|s-r\right|^{-2/3}\left(u_{r}\hat{P}(u_{r}),e_{-\lambda}\right)dr\right|^{p}\\
 & \leqslant C(\lambda,p,T)e_{\lambda p}(z)\bar{\delta}^{(1/3)(p-1)}\intop_{s-\bar{\delta}}^{s}\left|s-r\right|^{-2/3}\left(1,e_{-\lambda}\right)ds\\
 & \leqslant C(\lambda,p,T)\bar{\delta}^{p/3}e_{\lambda p}(z).
\end{align*}

Let us look at the 6th expectation:
\begin{align*}
 & E\left|\intop_{0}^{s-\bar{\delta}}(1,u_{r}^{2}(\psi_{t-r}^{z}-\psi_{s-r}^{z}))dr\right|^{p}\\
 & \leqslant E\left|\intop_{0}^{s-\bar{\delta}}\left(1,u_{r}^{2}e_{-\lambda}\right)dr\right|^{p}\sup_{r\in[0,s-\bar{\delta})}\left\Vert \psi_{t-r}^{z}-\psi_{s-r}^{z}\right\Vert _{\lambda}^{p}\\
 & \leqslant C(\lambda,p,T)e_{\lambda p}(z)(\left|t-s\right|^{p/2}\bar{\delta}^{-3p/2}+n^{-p/2}\bar{\delta}^{-3p/4}).
\end{align*}

Let us look at the 7th expectation:
\begin{align*}
 & E\left|\intop_{s-\bar{\delta}}^{s}(1,u_{r}^{2}(\psi_{t-r}^{z}-\psi_{s-r}^{z}))dr\right|^{p}\\
 & \leqslant C(p)E\left|\intop_{s-\bar{\delta}}^{s}\left\Vert \psi_{s-r}^{z}\right\Vert _{\lambda}\left(u_{r}^{2},e_{-\lambda}\right)ds\right|^{p}\\
 & \leqslant C(\lambda,p,T)e_{\lambda p}(z)\bar{\delta}^{(1/3)(p-1)}\intop_{s-\bar{\delta}}^{s}\left|t-s\right|^{-2/3}\left(1,e_{-\lambda}\right)ds\\
 & \leqslant C(\lambda,p,T)\bar{\delta}^{p/3}e_{\lambda p}(z).
\end{align*}

Let us look at the 8th expectation:
\begin{align*}
 & E\left|\intop_{0}^{s-\bar{\delta}}\left(v_{r},(\psi_{t-r}^{z}-\psi_{s-r}^{z})^{2}\right)dr\right|^{p/2}\\
 & \leqslant E\left|\intop_{0}^{s-\bar{\delta}}\left(v_{r},e_{-2\lambda}\right)dr\right|^{p/2}\sup_{r\in[0,s-\bar{\delta})}\left\Vert \psi_{t-r}^{z}-\psi_{s-r}^{z}\right\Vert _{\lambda}^{p}\\
 & \leqslant C(\lambda,p,T)e_{\lambda p}(z)(\left|t-s\right|^{p/2}\bar{\delta}^{-3p/2}+n^{-p/2}\bar{\delta}^{-3p/4}).
\end{align*}

Let us look at the 9th expectation:
\begin{align*}
 & E\left|\intop_{s-\bar{\delta}}^{s}\left(v_{r},(\psi_{t-r}^{z}-\psi_{s-r}^{z})^{2}\right)dr\right|^{p/2}\\
 & \leqslant E\left|\intop_{s-\bar{\delta}}^{s}\left\Vert \psi_{t-r}^{z}+\psi_{s-r}^{z}\right\Vert _{0}\left(v_{r},\psi_{t-r}^{z}+\psi_{s-r}^{z}\right)dr\right|^{p/2}\\
 & \leqslant C(T)E\left|\intop_{s-\bar{\delta}}^{s}\left|s-r\right|^{-2/3}\left(u_{r},\bar{\psi}_{t-r}^{z}+\bar{\psi}_{s-r}^{z}\right)dr\right|^{p/2}\\
 & \leqslant C(p,T)\bar{\delta}^{p/6}.
\end{align*}

Put all the above conclusions together:
\[
E(\left|\hat{u}_{t}(z)-\hat{u}_{s}(y)\right|^{p})\leqslant C(\lambda,p,T)e_{\lambda p}(z)(\left|t-s\right|^{p/24}+n^{-p/12}).
\]

\begin{lem}
For any $\lambda>0,T<\infty$
\end{lem}

\begin{enumerate}
\item $P\left(\underset{t\leqslant T}{\sup}\left\Vert \hat{u}_{t}(z)-\tilde{u}_{t}(z)\right\Vert _{-\lambda}\geqslant7n^{-1/4}\right)\rightarrow0$
as $n\rightarrow\infty$, 
\item $\underset{t\leqslant T}{\sup}\left\Vert (v_{0},\psi_{t}^{\cdot})-P_{t/3}f_{0}\right\Vert _{-\lambda}\rightarrow0$
as $n\rightarrow\infty$.
\end{enumerate}
The detail of proof is in \cite[Lemma 7]{mueller1995stochastic}.

From Kolmogorov's continuity criterion and the above moment estimate,
we can get the tightness of $\tilde{u}_{t}(z)$ as continuous $\mathscr{C}$
valued process. Then the tightness of $\hat{u}_{t}(z)$ follows, also
the continuity of all limit points follow from the above lemma.

\section{Characterizing limit points}

Taking a continuous function $\phi:\mathbb{R}\longrightarrow\mathbb{R}$
with compact support, we define

\[
E_{t}^{(6)}(\phi):=(v_{t}^{n},\phi)-(u_{n}(t),\phi),
\]

where $E(\underset{t\leqslant T}{\sup}\left|E_{t}^{(6)}(\phi)\right|)\leqslant C\left\Vert D(\phi,\mathcal{D}_{n})\right\Vert _{\lambda}$.
From the tightness of $u_{n}(t)$, we can get the tightness of $v_{t}^{n}$
as cadlag Radon measure valued process with the vague topology once
a compact containment condition is checked and all limit points are
again continuous.

Because of simultaneous convergence of $u_{n}(t)$ and $v_{t}^{n}$,
we write it in pair $(u_{n}(t),v_{t}^{n})$. By Skorokhod theorem,
we can find random variables with the same distribution as $\xi_{t}^{n}$,
which converges almost sure, and we still label it as $(u_{n}(t),v_{t}^{n})$.
Since the limits are continuous, the almost sure convergence holds
not only in Skorokhod sense but also in uniform sense on compact sets.
Thus, with probability one, for $T<\infty,\;\lambda>0,\phi$ of compact
support, we have:

\begin{align*}
 & \sup_{t\leqslant T}\left\Vert u_{n}(t)-u_{t}\right\Vert _{-\lambda}\longrightarrow0,\\
 & \sup_{t\leqslant T}\left\Vert \intop\phi(x)v_{t}^{n}(dx)-\intop\phi(x)v_{t}(dx)\right\Vert _{-\lambda}\longrightarrow0,
\end{align*}

where $v_{t}(dx)=u_{t}(x)dx$ for all $t\geqslant0$.

Pick up a $\phi$ three times continuously differentiable and with
compact support then substitute it for $\phi_{t}$ into (\ref{eq:1}):
we have:

\begin{align*}
Z_{t}(\phi) & =\intop\phi(x)v_{t}^{n}(dx)-\intop\phi(x)v_{0}^{n}(dx)\\
 & -\intop_{0}^{t}\intop(\Delta_{n}(\phi)(x)-\mathcal{L}_{n}^{c}\phi)v_{s}^{n}(dx)ds\\
 & -\intop_{0}^{t}\intop\left(-P(u_{s})\phi(x)+2\mathcal{L}_{n}^{c}\mathcal{S}_{n}u_{s}\phi(x)-2\mathcal{L}_{n}^{c}\mathcal{S}_{n}^{2}u_{s}^{2}\phi(x)\right)v_{s}^{n}(dx)ds\\
 & -\intop_{0}^{t}\intop\mathcal{L}_{n}^{c}\mathcal{S}_{n}^{2}u_{s}^{2}\phi(x)dxds-\sum_{i=1}^{5}E_{t}^{(i)}(\phi).
\end{align*}

When n goes to infinity,

$E_{t}^{(i)}(\phi)\;\left(1\leqslant i\leqslant5\right)$ tend to
zero for all t almost surely, from our assumption in case 1.

$\Delta_{n}(\phi)(x)$ tends to $\frac{1}{3}\Delta$ uniformly.

$\mathcal{L}_{n}^{c},\mathcal{S}_{n}$ are constants.

Therefore, $Z_{t}(\phi)$ tends to a continuous local martingale $z_{t}(\phi)$
where

\begin{align*}
z_{t}(\phi) & =\intop\phi(x)u_{t}(x)dx-\intop\phi(x)u_{0}(x)dx\\
 & -\intop_{0}^{t}\intop(\frac{1}{3}\Delta\phi(x)-2\phi(x)-P(u_{s})+4u_{s}(x)\phi(x)-4u_{s}^{2}(x)\phi(x)+2u_{s}(x)\phi(x))u_{s}(x)dxds.
\end{align*}

From (\ref{eq:2}), there exists 

\begin{align*}
 & Z_{t}^{2}(\phi)-2\intop_{0}^{t}\left(\mathcal{\rho}_{n}\mathcal{S}_{n}\right)^{-1}\mathcal{H}_{n}\left[(v_{s},\phi^{2})-(v_{s},\mathcal{S}_{n}\phi A(\xi_{t}^{n}\phi))\right]ds\\
= & Z_{t}^{2}(\phi)-2\intop_{0}^{t}\left(\mathcal{\rho}_{n}\mathcal{S}_{n}\right)^{-1}\mathcal{H}_{n}(v_{s},(1-\mathcal{S}_{n}u_{s})\phi^{2})ds\\
+ & 2\intop_{0}^{t}\left(\mathcal{\rho}_{n}\mathcal{S}_{n}\right)^{-1}\mathcal{H}_{n}(v_{s}\phi,\mathcal{S}_{n}A(\xi_{t}^{n}\phi_{s})-\mathcal{S}_{n}u_{s}\phi)ds\\
= & Z_{t}^{2}(\phi)-2\intop_{0}^{t}\left(\mathcal{\rho}_{n}\mathcal{S}_{n}\right)^{-1}\mathcal{H}_{n}(v_{s},(1-u_{s})\phi^{2})ds+E_{t}^{(7)}(\phi),
\end{align*}

where 
\[
E_{t}^{(7)}(\phi):=2\intop_{0}^{t}\left(\mathcal{\rho}_{n}\mathcal{S}_{n}\right)^{-1}\mathcal{H}_{n}(v_{s}\phi,\mathcal{S}_{n}A(\xi_{t}^{n}\phi)-\mathcal{S}_{n}u_{s}\phi)ds.
\]

Fortunately, we have a bound
\begin{align*}
E_{t}^{(7)}(\phi) & \leqslant2\intop_{0}^{t}\left(\mathcal{\rho}_{n}\mathcal{S}_{n}\right)^{-2}\mathcal{H}_{n}\sum_{x}\xi_{s}^{n}(x)\phi(x)(\mathcal{N}_{n}\mathcal{S}_{n})^{-1}\sum_{y\sim x}\mathcal{S}_{n}\xi_{s-}^{n}(y)\left|\phi(y)-\phi(x)\right|ds\\
 & \leqslant2t\left(\mathcal{\rho}_{n}\mathcal{S}_{n}\right)^{-1}\mathcal{H}_{n}(v_{t}\phi,\mathcal{S}_{n}u_{t}D(\phi,\mathcal{D}_{n}))\\
 & \leqslant2t\left(\mathcal{\rho}_{n}\mathcal{S}_{n}\right)^{-1}\mathcal{H}_{n}(\phi,D(\phi,\mathcal{D}_{n}))\rightarrow0.
\end{align*}

So
\[
z_{t}^{2}(\phi)-4\intop_{0}^{t}\intop(1-u_{s}(x))\phi^{2}(x)u_{s}(x)dxds
\]

is a continuous martingale.

Hence we have
\[
\partial_{t}u=\frac{1}{3}\text{\ensuremath{\Delta}}u+4u(u-\frac{1}{2})(1-u)-uP(u)+2\sqrt{u(1-u)}\dot{W.}
\]
\bibliographystyle{amsplain}
\bibliography{/Users/zhaotong/Downloads/paper/My_Bib}

\end{document}